\documentstyle{amsppt}


\magnification=1100
\vcorrection{-3truemm}
\pageheight{208truemm}



 

\author
{\smc PETER J. VASSILIOU}
\endauthor

\heading {Efficient Construction of Contact Coordinates for Partial Prolongations}\endheading
\vskip 20 pt
\heading{ Peter J. Vassiliou}\endheading
\centerline{ School of Mathematics and Statistics}
\centerline{ University of Canberra, Canberra, A.C.T., 2601}
\centerline{\tt peterv\@ise.canberra.edu.au}

\vskip 40 pt

\def\a{\alpha}
\def\b{\beta}
\def\d{\delta}
\def\D{\Delta}

\def\G{\Gamma}

\def\O{\Omega}
\def\f{\phi}

\def\o{\omega}
\def\O{\Omega}
\def\r{\rho}
\def\s{\sigma}
\def\S{\Sigma}

\def\CK{{\Cal K}}

\def\CV{{\Cal V}}
\def\P{\partial}

\def\P #1{\partial_{#1}}

\def\ch #1{\text{\rm Char}\ #1}

\document


\noindent {\eightpoint {\bf Abstract.} Let $\CV$ be a vector field distribution on manifold $M$. We give an efficient algorithm for the construction of local coordinates on $M$ such that $\CV$ may be locally expressed as some partial prolongation of the contact distribution 
$\Cal C^{(1)}_q$, on the $1^{st}$ order jet bundle of maps from $\Bbb R$ to $\Bbb R^q$, $q\geq 1$.   
It is proven that if $\CV$ is locally equivalent to a partial prolongation of $\Cal C^{(1)}_q$ then the explicit construction of contact coordinates algorithmically depends upon the determination of certain first integrals in a sequence of geometrically defined and algorithmically determined integrable Pfaffian systems on $M$. The number of these first integrals that must be computed satisfies a natural minimality criterion. These results therefore provide a full and constructive generalisation of the classical Goursat normal form from the theory of exterior differential systems.  

\vskip 20 pt
\item{}{\it 2000 MSC}: {\bf 58J60, 34H05, 93B18, 93B27.}
\item{}{\it keywords}: {Generalised Goursat normal form, partial prolongation, algorithm, nonlinear control theory, 
contact coordinates, Brunovsky normal form} }

\vskip 20 pt

\centerline{\bf 1. Introduction.}
\vskip 10 pt
Recent years have witnessed progress in the creation of geometric characterisations of contact distributions over jet bundles [1,9,4]. The problem of obtaining a simple, invariant characterisation of the corresponding {\it partial prolongations} is a natural one and has apparently remained open, even in the case of partial prolongations of jet bundles of maps from the line. That is, partial prolongations of the contact distribution $\Cal C^{(1)}_q$ over the first order jet bundle of maps $\Bbb R\to \Bbb R^q$, $q\geq 1$. 

In [6], a solution to the characterisation problem for partial prolongations of 
$\Cal C^{(1)}_q$ was presented. Specifically, simple geometric conditions, expressed in terms of the derived type of 
$\CV\subset TM$ were obtained, guaranteeing the existence of a local diffeomorphism from $M$ which identifies $\CV$ with some partial prolongation of $\Cal C^{(1)}_q$. 

However, for applications to integrable systems, nonlinear control theory and related areas, it is desirable not only to settle the recognition problem, but additionally, to find 
a method for explicitly {\it constructing} an equivalence between a differential system and a partial prolongation, whenever one is known to exist.

The main aim of this paper is to describe an efficient algorithm for this very construction problem.

We identify canonical and algorithmically determined integrable Pfaffian systems over the ambient manifold whose first integrals determine local equivalences. The number of these first integrals that must be constructed, according to our algorithm, satisfies a natural minimality criterion. In the special case when $\CV$ arises from a time-invariant nonlinear control system, the complexity of our algorithm agrees with that of the well-known GS algorithm [3]. However, even here, the method described in this paper enjoys a number of advantages over the former. This point is briefly discussed in section 4. We go on in section 5 to give an illustrative example of a nonlinear control system where, in addition, a new necessary condition for static feedback linearisation is derived.
 
Throughout the present work we will rely on the results of the aforementioned paper [6] and, wherever possible, adhere closely to its notational conventions. In sections 2 and 3, we will recall briefly the salient details of [6] that are required to establish the main result which is given in section 4. Section 5 is devoted to illustrative examples and an application.    

\vskip 25 pt

\centerline{\bf 2. Preliminaries}

 \vskip 5 pt
 {\it 2.1. The derived flag.}  The structure tensor of a sub-bundle $\CV\subset TM$ over manifold $M$ is a 
map
$$
\d_0\:\CV\times\CV\to TM/\CV,
$$
defined by
$$
\d_0(X,Y)={\bold p_0}([X,Y])
$$
where 
$$
{\bold p_0}\:TM\to TM/\CV
$$
is the natural projection.
If $\d_0$ has constant rank on $M$ we define the 
{\it first derived bundle} $\CV^{(1)}$ of $\CV$ as the maximal sub-bundle of $TM$ such that
$$
{\bold p_0}(\CV^{(1)})=\text{Im}~\d_0.
$$
In the constant rank case, we can iterate this procedure, obtaining the sequence of structure tensors
$$
\d_i\:\CV^{(i)}\times \CV^{(i)}\to TM/\CV^{(i)},\ i\geq 1
$$
each assumed to have constant rank on $M$. Define the sub-bundle $\CV^{(i)}\subset TM$,$\ i\geq 1$,
as the maximal sub-bundle of $TM$ satisfying ${\bold p_{i-1}}(\CV^{(i)})=\text{Im}~\d_{i-1}$
where ${\bold p_{i-1}}\:TM\allowbreak\to TM/\CV^{(i-1)}$ is the natural projection.
The sub-bundle $\CV^{(i)}\subset TM$ is the $i^{th}$ {\it derived bundle} of $\CV$. We call a bundle $\CV\subset TM$
{\it regular} if the ranks of all structure tensors are constant on $M$. In this regular case, for dimension reasons there is an integer $k\geq 0$ such that $\CV^{(k+1)}=\CV^{(k)}$. The smallest such integer is called the {\it derived length}
of $\CV$ and we evidently have a flag of sub-bundles
$$
\CV\subset\CV^{(1)}\subset\CV^{(2)}\subset\cdots\subset\CV^{(k)}
$$
called the {\it derived flag} of $\CV$. Plainly, the numbers $\dim\CV^{(i)}$ are diffeomorphism invariants of $\CV$.

\vskip 10 pt
{\it 2.2.  Cauchy bundles.} The Cauchy system or characteristic system $\chi(\CV)$ of $\CV$ is defined by
$$
\chi(\CV)=\big\{X\in\G(\CV)~|~[X,\G(\CV)]\subseteq\G(\CV)\big\}.
$$
Even if $\CV$ is regular, its Cauchy system need not have constant rank on $M$. However, if it does, there is a 
sub-bundle $\text{Char}~\CV\subseteq\CV$, the {\it Cauchy bundle}, whose space of smooth sections is $\chi(\CV)$.
\vskip 5 pt
\noindent{\smc Definition 2.1.} If a regular bundle $\CV\subset TM$ satisfies the condition that $\chi(\CV^{(i)})$ has constant rank on $M$ for each $i\geq 0$ then $\CV$ is said to be {\it totally regular}. For each $i$, the bundle 
$\text{Char}~\CV^{(i)}$ is called the {\it Cauchy bundle} or {\it characteristic bundle} of $\CV^{(i)}$.
The elements of $\text{Char}~\CV^{(i)}$ are called {\it Cauchy vectors}.
\vskip 5 pt
It is easy to show that any Cauchy bundle is integrable.
\vskip 5 pt
\noindent{\smc Definition 2.2.} Let $\CV\subset TM$ be a totally regular bundle. By the {\it derived type}
of $\CV$ we shall mean the list of sub-bundles \big\{$\CV^{(l)},\text{Char}~\CV^{(l)}\big\}$ for all
$l\geq 0$.
\vskip 10 pt
{\it 2.3. The singular variety of a sub-bundle.} The final tool we require is less well known than the derived type of a bundle, at least in the form in which we shall use it, so we discuss this in more detail in the remainder of this section. 

Let $\CV\subset TM$ be such that its structure tensor $\d_0$ has constant rank on $M$. For the remainder of this section we will not need to retain the subscript 0. For simplicity of notation, we will henceforth cease to distinguish between a bundle and its module of smooth sections. By $\Bbb P\CV$ we denote the projectivisation of $\CV$. 

Fix a basis $X_1,\ldots,X_n$ for $\CV$. If $Z_1,\ldots,Z_s$ complete $X_1,\ldots,X_n$
to a frame on $M$, then there are functions $c_{\a\b}^k$ on $M$, antisymmetric in lower indices, such that
$$
\d(X_\a,X_\b)=c_{\a\b}^k \bar{Z}_k, \ \ 1\leq \a,\b\leq n,\ 1\leq k\leq s,
$$
where $\bar{Z}_k$ is the coset with representative $Z_k$. The elements of $\Bbb P\CV$ are lines determined by elements 
$e^\a X_\a\in\CV$ and will be denoted by the symbol $E=[e^\a X_\a]$.
If $E$ is a fixed direction in $\CV$ then we may seek directions $F=[f^\a X_\a]$ such that
$$
\d(e^\a X_\a,f^\b X_\b)=0.\eqno(2.1)
$$
Equation (2.1) can be expressed in the form
$$
\sum_{\a<\b=2}^n(e^\a f^\b-e^\b f^\a)c^k_{\a\b}=0,\ 1\leq k\leq s,\eqno(2.2)
$$
viewed as a homogeneous linear system for $F$ and in the alternative form
$$
\s(E)\cdot{\bold f}=0,\eqno(2.3)
$$
where ${\bold f}=(f^1,f^2,\ldots,f^n)$. Letting ${\bold e}=(e^1,e^2,\ldots,e^n)$ we have clearly that 
${\bold f}={\bold e}$ is always a solution of (2.3).
\vskip 5 pt
{\smc Definition 2.3.} The matrix $\s(E)$ in (2.3) determined by equations (2.2) will be called the {\it polar matrix}
of the point $E\in\Bbb P\CV$. Such a line  will be called {\it singular} if its polar matrix has less than generic
rank. The set of all singular lines in $\CV$ will be denoted by the symbol $\text{Sing}(\CV)$. We will denote the set of singular lines in $\CV$ over a point $x\in M$, by $\text{Sing}(\CV)(x)$.
\vskip 5 pt
{\smc Lemma 2.1.} {\it For each $x\in M$, $\text{Sing}(\CV)(x)\subset \Bbb P\CV_x$, is a linear determinantal variety.}
\vskip 2 pt
\noindent {\it Proof.} On the one hand, for each $x\in M$, the polar matrix $\s(E_x)$ has entries which are linear functions
of the affine coordinates in $\Bbb P\CV_x$. On the other, $\text{Sing}(\CV)(x)$ is determined by equating the  minors of 
$\s(E_x)$ to zero.\hfill\hfill\qed
\vskip 5 pt
{\smc Definition 2.4.} The set $\text{Sing}(\CV)$ of all singular points of $\Bbb P\CV$, will be called the {\it singular variety} of $\CV$.
\vskip 5 pt
We now describe a very useful invariant associated to each point of ${\Bbb P}\CV$. 

The structure tensor $\d$  of $\CV\subset TM$ induces a map
$$
\text{deg}_{\CV}\:\Bbb P\CV\to\Bbb N
$$
well defined by
$$
\text{deg}_{\CV}([X])=\dim~\text{Image}~\D_X,\ [X]\in \Bbb P\CV
$$
where for each $X\in\CV$,
$$
\D_X\:\CV\to TM/\CV.
$$
is defined by 
$$
\D_X(Y)=\d(X,Y),\ \ Y\in\CV.
$$
\vskip 5 pt
{\smc Definition 2.5.} Let $\CV\subset TM$ be a sub-bundle and $[X]\in\Bbb P\CV$. The integer
$\text{deg}_\CV([X])$ will be called the {\it degree} of $[X]$. 
\vskip 5 pt
{\smc Lemma 2.2.} {\it For any $E\in\Bbb P\CV$, $\text{deg}_\CV(E)$ is a diffeomorphism invariant.}
\vskip 5 pt
Note that if $\text{deg}_\CV([X])=0$, then $X$ is a Cauchy vector and hence Lemma 2.2 implies, as a special case, the elementary fact that $\text{Char}~\CV$ is an invariant sub-bundle of $\CV$.

\vskip 5 pt
{\smc Lemma 2.3.} {\it For any point $E\in\Bbb P\CV$, $\text{deg}_\CV(E)=\text{rank}~\s(E)$. }
\vskip 2 pt
\noindent{\it Proof.} In fixed bases for $\CV$ and $TM/\CV$, the polar matrix $\s(E)$ is the matrix of the vector bundle 
morphism $\Delta_X$.\hfill\hfill\qed

\vskip 5 pt
{\smc Remark 2.1.} Amplifying Lemma 2.2, $\text{Sing}(\CV)$ is a diffeomorphism invariant in the sense that if $\CV_1,\CV_2$ are sub-bundles over
$M_1,M_2$, respectively and there is a diffeomorphism $\phi\: M_1\to M_2$ that identifies them, $\phi_*\CV_1=\CV_2$, then 
$\text{Sing}(\CV_2)$ and $\text{Sing}(\phi_*\CV_1)$ are equivalent as projective subvarieties of $\Bbb P\CV_2$. That is, for each $x\in M_1$, there is an element of the projective linear group $PGL({\CV_2}_{|_{\f(x)}},\Bbb R)$ that identifies 
$\text{Sing}(\CV_2)(\f(x))$ and $\text{Sing}(\phi_*\CV_1)(\f(x))$.
\vskip 5 pt
We hasten to point out that the computation of the singular variety for any given sub-bundle  $\CV\subset TM$ is algorithmic. That is, it involves only differentiation and commutative algebra operations. In practice, one computes the determinantal variety of the generic polar matrix $\s(E)$ as a sub-variety of $\Bbb P\CV$. 
\vskip 10 pt

{\it 2.4. The singular variety in positive degree.} One frequently encounters bundles possessing non-trivial Cauchy bundles and it is then natural and prudent to consider the notion of degree of a line in the 
quotient by $\ch\CV$. Continue with sub-bundle $\CV\subset TM$  and let
$$
\widehat{\pi}\:TM\to TM/\ch\CV:=\widehat{TM}
$$
be the natural projection assigning an element of $TM$ to its coset. The structure tensor $\d$ of $\CV$ drops to the quotient bundle $\widehat{\CV}:=\CV/\ch\CV$, leading to a tensor
$$
\widehat{\d}\:\widehat{\CV}\times\widehat{\CV}\to\widehat{TM}/\widehat{\CV}
$$   
well defined by
$$
\widehat{\d}(\widehat{X},\widehat{Y})=\widehat{\pi}\big([X,Y]\big)\mod\widehat{\CV},\ \ \forall\ 
\widehat{X},\widehat{Y}\in\widehat{\CV}
$$
where $\widehat{X}=\widehat{\pi}(X)$, $\widehat{Y}=\widehat{\pi}(Y)$. By analogy, we call $\widehat{\d}$ the {\it structure tensor} of $\widehat{\CV}$. 
As before, we are able to introduce a map
$$
\text{deg}_{\widehat{\CV}}\:\Bbb P\widehat{\CV}\to\Bbb N
$$
well defined by
$$
\text{deg}_{\widehat{\CV}}([\widehat{X}])=\dim~\text{Image}~\widehat{\D}_{\widehat{X}}
$$
where for $\widehat{X}\in\widehat{\CV}$,
$$
\widehat{\D}_{\widehat{X}}(\widehat{Y})=\widehat{\d}(\widehat{X},\widehat{Y}),\ \ \forall\ \ \ \widehat{Y}\in\widehat{\CV}.
$$
Here, by $[\widehat{X}]$ we denote the distribution spanned by $X\in\CV$ and $\ch\CV$. 
Again, by analogy we call $\text{deg}_{\widehat{\CV}}([\widehat{X}])$ the {\it degree} of 
$[\widehat{X}]\in\Bbb P\widehat{\CV}$. 
\vskip 5 pt
{\smc Remark 2.2.} All definitions and results of subsection 2.3 hold {\it mutatis mutandis} when the structure tensor $\d$  is replaced by
$\widehat{\d}$. In particular, we have notions of polar matrix and singular variety, as before. Note however, that each point of $\Bbb P\widehat{\CV}$ has degree one or more. 
\vskip 10 pt

{\it 2.5. The resolvent bundle.} In [6], we introduced the notion of a {\it Weber structure} which we now recall. 
Suppose a totally regular sub-bundle $\CV\subset TM$ of rank $c+q+1$, $q\geq 2, c\geq 0$ is defined on manifold $M$, 
$\dim M=c+2q+1$. Suppose further that $\CV$ satisfies the following additional properties:
\vskip 2 pt
{\itemitem{i)} $\dim\ch\CV=c$, $\CV^{(1)}=TM$,
\itemitem{ii)} 
$\widehat{\S}:=\text{Sing}(\widehat{\CV})=\Bbb P\widehat{\Cal B}\approx\Bbb R\Bbb P^{q-1}$ for some rank $q$ sub-bundle $\widehat{\Cal B}\subset\widehat{\CV}$}.
   
\vskip 5 pt
{\smc Definition 2.6.} We will call  $(\CV,\Bbb P\widehat{\Cal B})$ (or $(\CV,\widehat{\S})$) satisfying the above conditions a
{\it   Weber structure} on $M$. 
\vskip 5 pt         
Given a  Weber structure
$(\CV,\Bbb P\widehat{\Cal B})$, let $\Cal R_{\widehat{\S}}(\CV)\subset\CV$, denote the largest sub-bundle such that
$$
\widehat{\pi}\big( \Cal R_{\widehat{\S}}(\CV) \big)= \widehat{\Cal B}.\eqno(2.8)
$$
\vskip 5 pt
{\smc Definition 2.7.} The rank $q+c$ bundle $\Cal R_{\widehat{\S} }(\CV)$ will be called the {\it resolvent bundle} associated to the  Weber structure $({\CV},\widehat{\S})$. The bundle $\widehat{\Cal B}$ determined by the singular variety of 
$\widehat{\CV}$ will be called the {\it singular sub-bundle} of the  Weber structure. A  Weber structure
will be said to be {\it integrable} if its resolvent bundle is integrable.
\vskip 5 pt
{\smc Remark 2.3.} Note that an {\it integrable}  Weber structure descends to the quotient of $M$ by the leaves of $\ch\CV$ to be the contact bundle on $J^1(\Bbb R,\Bbb R^q)$. The term `Weber structure' honours Eduard von Weber (1870 - 1934) who, as far as I can tell, was the first to publish a proof of the Goursat normal form [7].
\vskip 5 pt
We record the following properties of the resolvent bundle of a Weber structure.
\vskip 5 pt
{\smc Proposition 2.7.}\ ([6])\ {\it  Let $(\CV,\widehat{\S})$ be a Weber structure on $M$ and 
$\widehat{\Cal B}$ its singular sub-bundle. If $q\geq 3$, then the following are equivalent
\vskip 3 pt
{\itemitem{1)} Its resolvent bundle $ \Cal R_{\widehat{\S}}(\CV) \subset\CV$ is integrable.                   
\itemitem{2)} Each point of $\widehat{\S}=\text{Sing}(\widehat{\CV})$ has degree one
\itemitem{3)} The structure tensor $\widehat{\d}$ of $\widehat{\CV}$ vanishes on $\widehat{\Cal B}$}}
\vskip 2 pt
{\smc Proposition 2.8.}\ ([6])\  {\it Let $({\CV},\widehat{\S})$ be an integrable Weber structure on $M$. 
Then its resolvent bundle $\Cal R_{\widehat{\S}}(\CV)$ is the unique, maximal, integrable sub-bundle of $\CV$.}

\vskip 5 pt
{\smc Remark 2.4.} Note that checking the integrability of the resolvent bundle is algorithmic.
One computes the singular variety $\text{Sing}(\widehat{\CV})=\Bbb P\widehat{\Cal B}$.
In turn, the singular bundle $\widehat{\Cal B}$ algorithmically determines  $\Cal R_{\widehat{\S}}(\CV)$.
\vskip 5 pt

We conclude this section by mentioning some notation. Firstly, in this paper we work exclusively in the smooth ($C^\infty$) category and all objects and maps will be assumed to be smooth without further notice. Secondly, we will often denote sub-bundles 
$\CV\subset TM$ by a list of  vector fields $X,Y,Z,\ldots$ on $M$ enclosed by braces,
$$
\CV=\big\{X,Y,Z,\ldots\big\}.
$$ 
This will always denote the bundle $\CV$ whose space of sections is the $C^\infty(M)$-module generated by
vector fields $X,Y,Z,\ldots$.

\vskip 25 pt
   
\centerline{\bf 3. Partial prolongations and Goursat bundles}
\vskip 10 pt

The {\it contact distribution} on the $1^{st}$ order jet bundle of maps from 
$\Bbb R\to \Bbb R^q$, $J^{1}(\Bbb R,\Bbb R^q)$, $q\geq 1$ will be denoted by the symbol $\Cal C^{(1)}_q$ and locally expressed in contact coordinates as
$$
\Cal C^{(1)}_q=\big\{\P x+\sum_{j=1}^qz^j_1\P {z^j}, \P {z^j_1}\big\}.\eqno(3.1)
$$
A {\it partial prolongation} of $\Cal C^{(1)}_q$ may be expressed in contact coordinates as a distribution on 
$J^k(\Bbb R,\Bbb R^q)$ of the form
$$
\Cal C(\tau)=\big\{\P x+\sum_{a=1}^t\sum_{j_a=1}^{q_a}\sum_{l_a=0}^{k_a-1} z^{j_a,a}_{l_a+1}\P {z^{j_a,a}_{l_a}}, 
\P {z^{j_a,a}_{k_a}}  \big\}.\eqno(3.2)
$$    
Here and elsewhere in this paper, the symbol $\tau$ denotes the {\it type} of the partial prolongation which may be specified by an ordered list of ordered pairs
$$
\tau=\langle q_1,k_1; q_2,k_2;\ldots;q_t,k_t\rangle\eqno(3.3)
$$
indicating that there are $q_a$ variables of order $k_a$ and we use the convention $1\leq k_1<k_2<\cdots <k_t$; the 
$q_a$ are any positive integers. The positive integer $t$ will be called the {\it class} of $\Cal C(\tau)$.

It will often be much more convenient to express the type of a partial prolongation as an ordered list of $k:=k_t$ non-negative integers 
$$
\tau=\langle \r_1,\r_2,\ldots,\r_k\rangle\eqno(3.4)
$$
where the $j^{th}$ element $\r_j$ indicates the number of variables of order $j$. In this notation, if $\r_j=0$ 
for all $j$ in the range 
$1\leq j\leq k-1$, then the contact distribution
(3.2) has class 1; its type is $\langle 0,0,\ldots,0,q\rangle$, $q\geq 1$. Such a contact system is a {\it total prolongation} of 
$\Cal C^{(1)}_q$, denoted, $\Cal C^{(k)}_q$, an {\it instance} of a partial prolongation; namely, the (standard) contact distribution on $J^k(\Bbb R,\Bbb R^q)$.  Note that, in the general case, the derived length of a partial prolongation of type $\tau$ is $k=k_t$.  

For any totally regular sub-bundle $\CV\subset TM$, we have the notion of its derived type. In section 2, we defined the {\it derived type} of a bundle as the list of all derived bundles together with their corresponding Cauchy bundles. We shall frequently abuse notation by using the term `derived type of $\CV$' for  the ordered list of ordered pairs of the form
$$
[[m_0,\chi^0],[m_1,\chi^1],\ldots,[m_k,\chi^k]]
$$
where $m_j$ denotes the rank of the $j^{th}$ derived bundle $\CV^{(j)}$ of $\CV$ and $\chi^j$ denotes the rank of its Cauchy bundle, $\ch\CV^{(j)}$. 

It is important to relate the type of a partial prolongation to its derived type. For this it's convenient to introduce the notions of {\it velocity}, {\it acceleration} and {\it decceleration} of a sub-bundle. 

\vskip 5 pt

{\smc Definition 3.1.} Let $\CV\subset TM$ be a totally regular sub-bundle of derived length $k$ and derived type 
$$
\d_\CV=[[m_0,\chi^0],[m_1,\chi^1],\ldots,[m_k,\chi^k]].\eqno(3.5)
$$ 
The {\it velocity} of $\CV$ is the ordered list of $k$ integers 
$$
\text{\rm vel}(\CV)=\langle\D_1,\D_2,\ldots,\D_k\rangle, 
$$
where,
$$
\D_j=m_j-m_{j-1},\ 1\leq j\leq k.
$$
The {\it acceleration} of $\CV$ is the ordered list of $k$ integers 
$$
\text{accel}(\CV)=\langle\D^2_2,\D^2_3,\ldots,\D^2_k,\D_k\rangle,
$$
where
$$
\D^2_i=\D_i-\D_{i-1},\ 2\leq i\leq k.
$$
The {\it decceleration} of $\CV$ is the ordered list of $k$ integers 
$$
\text{deccel}(\CV)=\langle -\D^2_2,-\D^2_3,\ldots,-\D^2_k, \D_k\rangle.
$$ 
Note that total prolongations $\Cal C^{(k)}_q$ have deccelerations of the form
$$
\langle 0,0,\ldots,0,q\rangle,\ q\geq 1,
$$
where there are $k-1$ zeros before the last entry $q$.  The classical Goursat normal form is the case $q=1$ in this family of deccelerations. The main result of [6] is a generalisation of this classical result to arbitrary deccelerations. The main aim of this paper is to describe an efficient algorithm for constructing the corresponding contact coordinates in this general case.
\vskip 5 pt
To recognise when a given sub-bundle has
or has not the derived type of a partial prolongation (3.2) we introduce one further canonically associated sub-bundle that plays a crucial role. 

If $\CV$ has derived length $k$ we let
$
\ch\CV^{(j)}_{j-1}
$
denote the intersections $\CV^{(j-1)}\cap\ch\CV^{(j)}$, $1\leq j\leq k-1$. It is easy to see that in every partial prolongation these sub-bundles are non-trivial and integrable.
\vskip 5 pt

{\smc Proposition 3.1.}\ ([6])\  {\it Let sub-bundle $\CV\subset TM$ be totally regular with velocity and acceleration
$\langle \D_1,\ldots,\D_k\rangle$ and $\langle\D^2_2,\ldots,\D^2_k,\D_k\rangle$, respectively. Then $\CV$ has the derived type
of a partial prolongation $\Cal C(\tau)$ of some type $\tau$ if and only if
$$
\aligned
&m_0=1+P,\  m_1=1+2P,\cr
&m_l=1+(1+l)P+\sum_{j=2}^l (l+1-j)\D^2_j,\ 2\leq l\leq k,\cr
&\chi^j=2m_j-m_{j+1}-1,\ 0\leq j\leq k-1,\cr
&\chi^i_{i-1}=m_{i-1}-1,\ 1\leq i\leq k-1.
\endaligned\eqno(3.6)
$$
where, $P=\sum_{i=1}^k\r_l$, and 
$$ 
m_l=\dim\CV^{(l)},\ \chi^l=\dim\ch\CV^{(l)},\ \chi^l_{l-1}=\dim\ch\CV^{(l)}_{l-1}. 
$$
The type $\tau$ in $\Cal C(\tau)$ is given by the decceleration,
$
\tau=\text{\rm deccel}(\CV).
$}
\vskip 5 pt

{\smc Definition 3.2.}  A totally regular sub-bundle $\CV\subset TM$ of derived length $k$ will be called a {\it Goursat bundle of type} $\tau$
if
\item{[i]} $\CV$ has the derived type of a partial prolongation whose type is $\tau=\text{deccel}(\CV)$ 
\item{[ii]} Each intersection $\ch\CV^{(i)}_{i-1}$ is an integrable sub-bundle whose rank, assumed to be constant on $M$,
agrees with the corresponding rank in $\Cal C(\tau)$
\item{[iii]} In case $\D_k>1$, then $\CV^{(k-1)}$ determines an integrable Weber structure on $M$, with singular sub-bundle
$\Bbb P\widehat{\Cal B}\approx\Bbb R\Bbb P^{\D_k-1}$.

\vskip 25 pt

\centerline{\bf 4. Efficient construction of contact coordinates}
\vskip 10 pt
In this section we give our main result, an efficient algorithm for the construction of contact coordinates for any smooth distribution locally equivalent to a partial prolongation of the contact distribution 
$\Cal C^{(1)}_q$ on $J^1(\Bbb R, \Bbb R^q)$, $q\geq 1$. Specifically, we show how to algorithmically construct certain canonical Pfaffian systems over the ambient manifold and appropriate first integrals of these Pfaffian systems whose derivatives generate a local equivalence, identifying a given Goursat bundle with some partial prolongation. The type of $\CV$ is given by the decceleration vector, deccel($\CV$). 

The main result upon which our algorithm is based is the generalised Goursat normal form established in [6].
\vskip 5 pt
{\smc Theorem 4.1.} [Generalised Goursat Normal Form], [6]. {\it Let $\CV\subset TM$ be a Goursat bundle over manifold $M$, of derived length $k>1$ and type 
$\tau=\text{\rm deccel}(\CV)$. 
Then there is an open, dense subset $\hat M\subseteq M$ such that the restriction of $\CV$ to $\hat M$ is locally equivalent to $\Cal C(\tau)$. Conversely any partial prolongation of $\Cal C^{(1)}_q$ is a Goursat bundle.}
\vskip 5 pt
This theorem settles the recognition problem for a distribution in terms of simple constraints on its derived type. 
While the  proof in [6] is constructive, it is extravagant with respect to the number of integrations that are carried out. The aim here is to show that the number of integrations that must be carried out to actually {\it construct} an equivalence for any Goursat bundle is comparatively small. In fact we show that the number of first integrals that must be computed in the general case is equal to the number of ``dependent variables" featured in $\Cal C(\tau)$, plus one, where
$\tau=\langle \rho_1,\rho_2,\ldots,\rho_k\rangle$ is the bundle's decceleration vector deccel($\CV$). That is,
$\sum_{j=1}^k\r_j+1=P+1$ first integrals must be found in the general case. The remaining coordinates are computed by differentiation. This is the ``natural" minimality criterion alluded to in the Introduction. Our algorithm therefore performs as well as the GS algorithm [3] in the special case when distribution $\CV$ happens to arise from a time-invariant control system. 

However, our approach is not restricted to time-invariant control systems, nor, indeed to general control systems. Furthermore, it doesn't involve the construction of structure equations
before an equivalence can be found. Moreover, the geometric data that must be computed to settle the recognition problem is expressed more naturally in terms of the bundle's derived type. We now procede to the description of our method.

\vskip 5 pt

Let $\CV\subset TM$ be a Goursat bundle over $M$ of derived length $k$. Recall that there is a distinction between the cases $\r_k>1$ and $\r_k=1$. In the former case we have the filtration 

$$
\aligned
\ch\CV^{(1)}_0&\subseteq\ch\CV^{(1)}\subset\cdots\subset\ch\CV^{(j)}_{j-1}\subseteq\ch\CV^{(j)}\subset\cdots\cr
 \cdots&\subset\ch\CV^{(k-1)}_{k-2}\subseteq\ch\CV^{(k-1)}\subset\Cal R_{\widehat{\S}_{k-1}}(\CV^{(k-1)})\subset TM
\endaligned\eqno(4.1)
$$
where $\Cal R_{\widehat{\S}_{k-1}}(\CV^{(k-1)})$ is the resolvent bundle of the integrable Weber structure 
$(\CV^{(k-1)},\allowbreak \widehat{\S}_{k-1})$. Note that for $j$ in the range $1\leq j\leq k-1$, $\ch\CV^{(j)}_{j-1}=\ch\CV^{(j)}$ if and only if $\D^2_{j+1}=0$. 

We will use the convenient notation
$$
\aligned
&\nu_j=\D_j,\ 1\leq j\leq k,\cr
&n_i=\nu_{i+1},\ 0\leq i\leq k-1,\cr  
&N_l=\dim M-m_l,\ 0\leq l\leq k.
\endaligned\eqno(4.2)
$$
Also, let
$$
\Xi^{(j)}={\ch\CV^{(j)}}^\perp,\ \Xi^{(j)}_{j-1}={\ch\CV^{(j)}_{j-1}}^\perp,\ 1\leq j\leq k-1
$$
and
$$
\Upsilon_{\widehat{\S}_{k-1}}(\CV^{(k-1)})={\Cal R_{\widehat{\S}_{k-1}}(\CV^{(k-1)})}^\perp.
$$
Then we have a filtration of the cotangent bundle $T^* M$
$$
\Upsilon_{\widehat{\S}_{k-1}}(\CV^{(k-1)})\subset\Xi^{(k-1)}\subseteq\Xi^{(k-1)}_{k-2}\subset\cdots\subset
\Xi^{(1)}\subseteq\Xi^{(1)}_0\subset T^*M\eqno(4.3)
$$
by integrable sub-bundles.
It follows easily from Proposition 3.1 and (4.2) that
$$
\aligned
&\dim\Upsilon_{\widehat{\S}_{k-1}}(\CV^{(k-1)})=N_{k-1}+1,\cr    
&\dim\Xi^{(j)}=N_j+\Delta_{j+1}+1,\ \dim\Xi^{(j)}_{j-1}=N_{j-1}+1,\ 1\leq j\leq k-1,
\endaligned\eqno(4.4)
$$
and 
$$
\dim\Xi^{(j)}_{j-1}-\dim\Xi^{(j)}=-\Delta^2_{j+1}=\rho_j,\ 1\leq j\leq k-1.\eqno(4.5)
$$
We can therefore construct a filtered basis for sub-bundle $\Xi^{(1)}_0\subset T^*M$ as follows
$$
\aligned
&\o_0,\o_1,\ldots,\o_{N_{j-1}}\ \ \text{for}\ \ \Xi^{(j)}_{j-1},\cr
&\o_0,\o_1,\ldots,\o_{N_j+n_j}\ \ \text{for}\ \ \Xi^{(j)},\cr
&\o_0,\o_1,\ldots,\o_{N_{k-1}}\ \ \text{for}\ \ \Upsilon_{\widehat{\S}_{k-1}}(\CV^{(k-1)}),
\endaligned\ \ 1\leq j\leq k-1.\eqno(4.6)
$$
{\smc Definition 4.1.} For each $j\in\{1,\ldots,k-1\}$ such that $\r_j>0$, we refer to the sub-bundle 
$\O_j=\{\o_{N_j+n_j+1},\ldots,\o_{N_j+\nu_j}\}\subset\Xi^{(j)}_{j-1}$ 
as the {\it fundamental bundle of order $j$}. 
\vskip 5 pt
For each $j$, for which $\r_j>0$, $\O_j$ is a sub-bundle of the integrable bundle $\Xi^{(j)}_{j-1}$ and we note the decomposition 
$\Xi^{(j)}_{j-1}=\Xi^{(j)}\oplus\O_j$. By the Frobenius theorem, there are functions
$\{\varphi^{l_j,j}\}_{l_j=1}^{\r_j}$ and multiplier matrices $M_j$ with $\det M_j\neq 0$ such that
$$
\left(\matrix d\varphi^{1,j}\\ d\varphi^{2,j}\\ \cdot \\ \cdot \\ d\varphi^{\rho_j,j}
\endmatrix\right)\equiv M_j\left(\matrix \o_{N_j+n_j+1} \\ \o_{N_j+n_j+2} \\ \cdot\\ \cdot \\ \o_{N_j+\nu_j}
\endmatrix\right)\mod\Xi^{(j)}.\eqno(4.7)
$$
\vskip 5 pt
{\smc Definition 4.2.} We refer to the functions $\varphi^{l_j,j}$, $j\in\{1,2,\ldots,k-1\}$, defined by (4.7) as  
{\it fundamental functions of order $j$}.
Let $\varphi^{0,k},\varphi^{1,k},\ldots,\varphi^{N_{k-1},k}$ span the first integrals of the integrable sub-bundle 
$\Upsilon_{\widehat{\S}_{k-1}}(\CV^{(k-1)})$. We refer to these as {\it fundamental functions of order $k$}. 
\vskip 5 pt
We now prove that the construction of the fundamental functions of all orders is the only integration that need be carried out in order to construct contact corrdinates for any Goursat bundle $\CV$.
\vskip 5 pt
{\smc Theorem 4.2.} {\it Let $\CV\subset TM$ be a Goursat bundle of derived length $k$ with decceleration vector 
$\tau=\text{\rm deccel}(\CV)=$ $\langle \rho_1,\rho_2,\ldots,\rho_k\rangle$, $\rho_k\geq 2$. 
Let $\{x,\varphi^{1,k},\varphi^{2,k},\ldots,\varphi^{\rho_k,k}\}$ denote the fundamental functions of order $k$
and for each $j$ in the range $1\leq j\leq k-1$ for which $\rho_j>0$, let $\{\varphi^{1,j},\ldots,\varphi^{\rho_j,j}\}$ denote the fundamental functions of order $j$ defined on some open subset $\Cal U\subseteq M$. 

Then there is an open, dense subset $\widehat{\Cal U}\subseteq \Cal U$ and a section $Y$ of $\CV$ such that on 
$\widehat{\Cal U}$, $Yx\neq 0$ and the fundamental functions $x,\varphi^{l_j,j}_0:=\varphi^{l_j,j}$, together with the functions
$$
\varphi^{l_j,j}_{s_j+1}=\frac{Y\varphi^{l_j,j}_{s_j}}{Yx},
\ j\in\{1,\ldots,k\},\ 1\leq l_j\leq \rho_j,\ 0 \leq s_j\leq j-1\eqno(4.8)
$$ 
are contact coordinates for $\CV$, identifying it with the partial prolongation $\Cal C(\tau)$.
 }

\vskip 3 pt

{\it Proof.} Fix a point $\bar{y}_0\in M$ in a neighbourhood $\Cal U$ of which $\CV$ is a Goursat bundle. The proof of Theorem 4.1 in [6] shows that we may extend the fundamental functions of order $j\in\{1,\ldots,k\}$, for which $\r_j>0$, 
namely, $x=\varphi^{0,k},z^{l_j,j}=\varphi^{l_j,j}$, $1\leq l_j\leq \r_j$, to a system of contact coordinates
$$
\bar{z}=\Big(x,z^{l_j,j},z^{l_j,j}_1,z^{l_j,j}_2,\ldots,
z^{l_j,j}_j\Big)_{j=1,l_j=1}^{k,\r_j}\eqno(4.9)
$$
on an open set $\widehat{\Cal U}\subseteq\Cal U$. Since $\ch\CV^{(1)}_0$ has codimension 1 in $\CV$, it follows that
there is a section $Z$ of $\CV$ such that $Zx=1$ on a dense open subset of $\widehat{\Cal U}$, which need not, in fact,  contain $\bar{y}_0$, and which, for simplicity, we denote by the same symbol.  Let
$$
\psi \: \bar{y} \mapsto \bar{z}
$$
be the local diffeomorphism defined by the change of variable from the original coordinates $\bar{y}$ to the contact coordinates (4.9). Then we have that $\psi_*\CV=\Cal C(\tau)$, where 
$$
\Cal C(\tau)=\{X=\P x+\sum_{j\in\{1,\ldots,k\}}\sum_{l_j=1}^{\r_j}\sum_{h_j=0}^{j-1} 
z_{h_j+1}^{l_j,j}\P {z^{l_j,j}_{h_j}},\P {z^{l_j,j}_j}\},
$$
and $z^{l_j,j}_0:=z^{l_j,j}$. Consequently, we have for some functions $\alpha,\alpha^{l_j,j}$ on $\widehat{\Cal U}$
$$
\psi_*Z=\alpha X+\sum_{l_j=1}^{\r_j}\alpha^{l_j,j}\P {z^{l_j,j}_{j}}
$$ 
Since
$$
\a=(\psi_*Z)x=Z(x\circ\psi)=Z\varphi^{0,k}=1
$$
we have
$$
\psi_*Z=X+\sum_{l_j=1}^{\r_j}\alpha^{l_j,j}\P {z^{l_j,j}_{j}}.\eqno(4.10)
$$ 
Define functions $\varphi^{l_j,j}_1,\varphi^{l_j,j}_2,\ldots,\varphi^{l_j,j}_j$ for  $j\in\{1,2,\ldots,k\}$
such that $\r_j>0$ by
$$
x=\varphi^{0,k},\varphi^{l_j,j}_{s_j+1}=Z\varphi^{l_j,j}_{s_j},\ 1\leq l_j\leq \r_j,\ 0\leq s_j\leq j-1,\eqno(4.11)
$$
where $\varphi^{l_j,j}_0:=\varphi^{l_j,j}$.
By (4.10) we have
$$
\varphi^{l_j,j}_1=Z(\varphi^{l_j,j})(\bar{y})=Z(\psi^*z^{l_j,j})(\bar y)=(\psi_*Z)(z^{l_j,j})(\psi(\bar y))=
\psi^*(z^{l_j,j}_1).\eqno(4.12)
$$
The calculation in (4.12) can be repeated so that for each $r_j$ in the range $1\leq r_j\leq j$, we have
$\varphi^{l_j,j}_{r_j}=\psi^*(z^{l_j,j}_{r_j})$ showing that the functions defined in (4.11) are independent
on $\widehat{\Cal U}$. By their definition, functions (4.11), together with $\varphi^{0,k},\varphi^{l_j,j}$,  are contact coordinates there. The proof is completed by observing that if $Y$ is any section of $\CV$ such that $Yx\neq 0$ on $\widehat{\Cal U}$ then we may take $Z$ to be the vector field
$(Yx)^{-1}Y$. \hfill\hfill\qed
\vskip 5 pt
The only case remaining is $\r_k=1$. Here  
$
(M/\ch\CV^{(k-1)},\CV/\ch\CV^{(k-1)})
$ 
is a 3-dimensional contact manifold so there is no canonical maximal integrable sub-bundle of 
$\CV^{(k-1)}$ as there is in the case $\r_k\geq 2$. The role of the resolvent bundle when $\r_k=1$ is played by a locally defined bundle, $\Pi^k$, whose construction we now describe.  Let $x$ denote any first integral of $\ch\CV^{(k-1)}$ and seek any section $Z$ of $\CV$ such that $Zx=1$. Define a sub-bundle 
$\Pi^k\subset\CV^{(k-1)}$ inductively by
$$
\Pi^{l+1}=[\Pi^l,Z],\ \Pi^0=\ch\CV^{(1)}_0,\ 0\leq l\leq k-1.\eqno(4.13)
$$
The proof of Theorem 4.1 shows that $\Pi^k$ is integrable, has codimension 2 in $TM$ and first integral $x$. 
In this case, filtration (4.3) is replaced by
$$
{\Pi^k}^\perp\subset\Xi^{(k-1)}\subseteq\Xi^{(k-1)}_{k-2}\subset\cdots\subset
\Xi^{(1)}\subseteq\Xi^{(1)}_0\subset T^*M\eqno(4.14)
$$
and the filtered basis (4.6) for $\Xi^{(1)}_0$ is replaced by 
$$
\aligned
&\o_0,\o_1,\ldots,\o_{N_{j-1}}\ \ \text{for}\ \ \Xi^{(j)}_{j-1},\cr
&\o_0,\o_1,\ldots,\o_{N_j+n_j}\ \ \text{for}\ \ \Xi^{(j)},\cr
&\o_0,\o_1,\ldots,\o_{N_{k-1}}\ \ \text{for}\ \ {\Pi^k}^\perp.
\endaligned\ \ 1\leq j\leq k-1,\eqno(4.15)
$$
Then by an argument similar to that of Theorem 4.2, we have
\vskip 5 pt
{\smc Theorem 4.3.} {\it Let $\CV\subset TM$ be a Goursat bundle of derived length $k$ and decceleration vector 
$\tau=\langle \r_1,\r_2,\ldots,\r_k\rangle$, $\r_k=1$. Let $\Pi^k$ be the bundle locally defined in (4.13).
Let $\varphi^{1,k}$ be any other first integral of $\Pi^k$
such that $dx\wedge d\varphi^{1,k}\neq 0$ on an open set $\Cal U\subseteq M$.

Then there is an open, dense subset 
$\widehat{\Cal U}\subseteq\Cal U$ upon which is defined a section $Z$ of $\CV$ satisfying $Zx=1$ such that the fundamental functions $x,\varphi^{l_j,j}_0:=\varphi^{l_j,j}$ together with the functions
$$
\varphi^{l_j,j}_{s_j+1}=X\varphi^{l_j,j}_{s_j},\ j\in\{1,\ldots,k\},\ 0\leq s_j\leq j-1,\ 1\leq l_j\leq\r_j,
$$
are contact coordinates for $\CV$ on $\widehat{\Cal U}$, indentifying it with the partial prolongation 
$\Cal C(\tau)$.}

\vskip 5 pt
Theorems 4.2 and 4.3 prove the correctness of algorithm {\tt Contact}. 

\newpage

\centerline{\bf Algorithm Contact A}
{\eightpoint
\centerline{------------------------------------------------------------------------------}
\vskip 5 pt
\item{} {\tt INPUT:} Goursat bundle $\CV\subset TM$ of derivedlength $k$ and type 
$\tau=\text{\rm deccel}(\CV)=\langle \r_1,\ldots,\r_k\rangle$, $\r_k>1$.
\item{1.} Build filtration (4.1) of $TM$ and filtration (4.3) of $T^*M$. 
\item{2.} Build filtered basis (4.6) of bundle $\Xi^{(1)}_0$ constructed in {\it step 1}.
\item{3.} For each $j$, if $\r_j>0$, compute the fundamental bundle $\O_j$ of order $j$.
\item{4.} For each $j$, if $\O_j$ is not empty, compute the fundamental functions $\{\varphi^{l_j,j}\}_{l_j=1}^{\r_j}$ of order $j$.
\item{5.} Fix any fundamental function of order $k$, denoted $x$ and any section $Z$ of $\CV$ such that $Zx=1$.
\item{6.} For each $j$, if $\r_j>0$ let $z^{l_j,j}=\varphi^{l_j,j}$, $1\leq l_j\leq \r_j$.
\item{7.} For each $j$, if $\r_j>0$ define functions
$$
x, z^{l_j,j}_0:=z^{l_j,j}=\varphi^{l_j,j},\ z^{l_j,j}_{s_j+1}=Zz^{l_j,j}_{s_j},\ 0\leq s_j\leq j-1,\ 1\leq l_j\leq\r_j.
$$
\item{} {\tt OUTPUT:} Contact coordinates for $\CV$ identifying it with $\Cal C(\tau)$.

\vskip 20 pt

\centerline{\bf Algorithm Contact B}

\centerline{-------------------------------------------------------------------------------}
\vskip 5 pt
\item{} {\tt INPUT:} Goursat bundle $\CV\subset TM$ of derived length $k$ and type 
$\tau=\text{\rm deccel}(\CV)=\langle \r_1,\ldots,\r_k\rangle$, $\r_k=1$.
\item{1.} Compute filtration (4.1) of $TM$ up to $\ch\CV^{(k-1)}$. 
\item{2.} Fix any first integral of $\ch\CV^{(k-1)}$, denoted $x$, and any section $Z$ of $\CV$ such that $Zx=1$ .
\item{3.} Build distribution $\Pi^k$, defined by (4.13), giving refinement (4.14) of the filtration constructed in 
{\it step 1}.
\item{4.} Let $z^k:=\varphi^{1,k}$ be any first integral of $\Pi^k$ such that $dx\wedge d\varphi^{1,k}\neq 0$. 
\item{5.} Build filtration (4.14) of $T^*M$ from that of $TM$ constructed in {\it steps 1 and 3}.
\item{6.} Build filtered basis (4.15) of bundle $\Xi^{(1)}_0$. 
\item{7.} For each $j$, if $\r_j>0$, compute the fundamental bundle $\O_j$ of order $j$.
\item{8.} For each $j$, if $\O_j$ is not empty, compute the fundamental functions $\{\varphi^{l_j,j}\}_{l_j=1}^{\r_j}$ of order $j$.
\item{9.} For each $j$, if $\r_j>0$ let $z^{l_j,j}=\varphi^{l_j,j}$, $1\leq l_j\leq \r_j$.
\item{10.} For each $j$, if $\r_j>0$ define functions
$$
x, z^{l_j,j}_0:=z^{l_j,j}=\varphi^{l_j,j},\ z^{l_j,j}_{s_j+1}=Zz^{l_j,j}_{s_j},\ 0\leq s_j\leq j-1,\ 1\leq l_j\leq \r_j.
$$
\item{} {\tt OUTPUT:} Contact coordinates for $\CV$ identifying it with $\Cal C(\tau)$.
\vskip 5 pt

}

\vskip 10 pt

\vskip 25 pt

\heading{\bf 5. Examples and Applications}\endheading
\vskip 5 pt
In this section we illustrate our algorithm for finding contact coordinates for distributions locally equivalent to partial prolongations and indicate some applications to nonlinear control theory. 
\vskip 5 pt
{\it Example 5.1.} We begin with a simple illustrative example that is sufficiently non-trivial to illustrate the all the steps in algorithm {\tt Contact A}. Consider the sub-bundle of $T\Bbb R^{21}$ defined by
{\eightpoint 
$$
\aligned
&\CV=\Big\{e^{-x_1}\P {x_1}+(x_4-x_8-e^{x_1})\P {x_2}+(x_5-x_{21}-x_{20}+x_4-x_8-e^{x_1})\P {x_3}+\cr
&(x_9+2x_{21}+2x_{20}-4x_4+4x_8+3e^{x_1})(\P {x_4}+\P {x_8})+(2-2x_{12}+x_7)\P {x_6}+\cr
&(x_8+e^{x_1})\P {x_7}-\P {x_8}+\P {x_9}+(1-x_{13}-x_{14}-x_{12})\P {x_{10}}+\cr
&(\frac{1}{2}(x_{13}+x_{14})+x_4-x_8-e^{x_1})\P {x_{11}}+\frac{1}{2}(x_{13}-x_{14})\P {x_{12}}+\cr
&(x_{15}-x_{16})\P {x_{13}}+(x_{15}+x_{16})\P {x_{14}}+x_{17}\P {x_{15}}+(x_{18}-2x_2)\P {x_{16}}+\cr
&(x_{19}-x_8-e^{x_1})\P {x_{17}}+x_{20}\P {x_{18}}+(3x_{21}+3x_{20}-6x_4+6x_8+5e^{x_1}+x_9)\P {x_{19}},\cr
&\P {x_4}+2\P {x_{20}}, \P {x_5},\P {x_9}, \P {x_{20}}-\P {x_{21}}, \P {x_5}-2\P {x_9}+\P {x_{21}}\Big\}
=\{X_1,X_2,\ldots,X_6\}.
\endaligned
$$
}
We compute the (refined) derived type to be
$$
[[6,0],[11,5,7],[14,10,10],[17,13,14],[19,16,16],[21,21]]
$$
where for $i=0$ and $i=5$ the two-element lists record $[m_i,\chi^i]$ and for $1\leq i\leq 4$ the 3-element lists
record $[m_i,\chi^i_{i-1},\chi^i]$. Hence the derived length is $k=5$ and by Proposition 3.1, we confirm that $\CV$ has the derived type of a partial prolongation whose type is $\langle 2,0,1,0,2\rangle$. Since $\r_5=2>1$, we check the singular variety of $\widehat{\CV}^{(4)}=\CV^{(4)}/\ch\CV^{(4)}$. We compute that
$$
\widehat{\CV}^{(4)}=\{\widehat{X}_1=\widehat{\pi}(X_1),\ \widehat{X}_2=\widehat{\pi}(\P {x_{12}}),\
\widehat{X}_3=\widehat{\pi}(\P {x_{13}}+\P {x_{14}})\}
$$
whose non-zero structure is
$$
\widehat{\d}(\widehat{X}_1,\widehat{X}_2)=\big\langle\widehat{\pi}(2\P {x_6}+\P {x_{10}})\big\rangle,\ 
\widehat{\d}(\widehat{X}_1,\widehat{X}_3)=\big\langle\widehat{\pi}(2\P {x_{10}}-\P {x_{11}})\big\rangle.
$$
Here $\widehat{\d}$ is the structure tensor of $\widehat{\CV}^{(4)}$,
$$
\widehat{\pi} : T\Bbb R^{21}\to T\Bbb R^{21}/\ch\CV^{(4)}=:\widehat{T}\Bbb R^{21}
$$
is the natural projection and for $Z\in T\Bbb R^{21}$, $\big\langle\widehat{\pi}(Z)\big\rangle$ denotes the element of the quotient bundle
$$
\widehat{T}\Bbb R^{21}/\widehat{\CV}^{(4)}
$$
with representative $\widehat{\pi}(Z)$.
From this, we 
easily compute that an arbitrary point $E=[a^1\widehat{X}_1+a^2\widehat{X}_2+a^3\widehat{X}_3]\in\Bbb P\widehat{\CV}^{(4)}$ has polar matrix
$$
\left(\matrix -a^2 & a^1 & 0\\
              -a^3 & 0 & a^1
      \endmatrix\right),
$$
from which we deduce that
$$
\text{\rm Sing}(\widehat{\CV}^{(4)})=\Bbb P\{\widehat{\pi}(\P {x_{12}}),\widehat{\pi}(\P {x_{13}}+\P {x_{14}})\}
$$
and the resolvent bundle is therefore
$$
\Cal R_{\widehat{\S}_{4}}(\CV^{(4)})=\ch\CV^{(4)}\oplus\{\P {x_{12}},\P {x_{13}}+\P {x_{14}}\}.
$$
The filtration induced on $T\Bbb R^{21}$ by $\CV$ is
$$
\aligned
\ch\CV^{(1)}_0\subset\ch\CV^{(1)}\subset\ch\CV^{(2)}\subset&\ch\CV^{(3)}_2\subset\ch\CV^{(3)}\cr
                                       &\subset\ch\CV^{(4)}\subset\Cal R_{\widehat{\S}_{4}}(\CV^{(4)})\subset T\Bbb R^{21}
\endaligned\eqno(5.1)
$$
and one can check that the resolvent bundle is integrable. This data confirms that $\CV$ is a Goursat bundle of type
$\text{\rm deccel}(\CV)=$ $\langle 2,0,1,0,2\rangle$.  By Theorem 4.1, $\CV$ is locally equivalent to the partial prolongation with this type. We now go on to use algorithm {\tt Contact} to construct an explicit equivalence.   

The filtration of $T^*\Bbb R^{21}$ induced by $\CV$ and dual to (5.1) is, in this case,
$$
\Upsilon_{\widehat{\S}_4}(\CV^{(4)})\subset\Xi^{(4)}\subset\Xi^{(3)}\subset\Xi^{(3)}_2\subset\Xi^{(2)}
                                       \subset\Xi^{(1)}\subset\Xi^{(1)}_0\subset T^*\Bbb R^{21}
\eqno(5.2)
$$
Pursuing {\it step 2} of algorithm {\tt Contact A}, a filtered basis for $\Xi^{(1)}_0$ is 
$$
\aligned
&\Upsilon_{\widehat{\S}_4}(\CV^{(4)})=\{dx_1,dx_2-dx_{11},dx_{10}\}\cr
&\Xi^{(4)}=\{dx_1,dx_2-dx_{11},dx_{10},dx_{12},dx_{13}+dx_{14}\}\cr
&\Xi^{(3)}=\{dx_1,dx_2-dx_{11},dx_{10},dx_{12},dx_{13}+dx_{14},dx_{15},dx_{14}\}\cr
&\Xi^{(3)}_2=\{dx_1,dx_2-dx_{11},dx_{10},dx_{12},dx_{13}+dx_{14},dx_{15},dx_{14},dx_6\}\cr
&\Xi^{(2)}=\{dx_1,dx_2-dx_{11},dx_{10},dx_{12},dx_{13}+dx_{14},dx_{15},dx_{14},dx_6,dx_7,dx_{16},dx_{17}\}\cr
&\Xi^{(1)}=\{dx_1,dx_2-dx_{11},dx_{10},dx_{12},dx_{13}+dx_{14},dx_{15},dx_{14},dx_6,dx_7,dx_{16},dx_{17},\cr
& \hskip 250 pt  dx_{19},dx_8,2dx_{11}-dx_{18}\}\cr
&\Xi^{(1)}_0=\{dx_1,dx_2-dx_{11},dx_{10},dx_{12},dx_{13}+dx_{14},dx_{15},dx_{14},dx_6,dx_7,dx_{16},dx_{17},\cr
& \hskip 220 pt  dx_{19},dx_8,2dx_{11}-dx_{18},dx_3,dx_{18}\}
\endaligned
$$
Since $\r_1,\r_3$ and $\r_5$ alone are non-zero, we deduce from this that the fundamental bundles are
$$
\O_1=\{dx_3,dx_{18}\},\ \O_3=\{dx_6\},\ \O_5=\{dx_1,dx_2-dx_{11},dx_{10}\}
$$
and the rest are empty. The corresponding fundamental functions may therefore be taken to be (for instance)
$$
\Cal F_1=\{z^{1,1}=x_3,z^{2,1}=x_{18}\},\ \Cal F_3=\{z^{1,3}=x_6\},\ \Cal F_5=\{z^{1,5}=-x_2+x_{11},z^{2,5}=x_{10}\}
$$
and we take the independent variable $x$ to be $e^{x_1}$, for then $X_1e^{x_1}=1$. Consequently, $Z=X_1$ may be taken to be the operator of total differentiation.

Finally, executing {\it step 7}, we compute the remaining components of the equivalence by differentiating once the elements of 
$\Cal F_1$, differentiating three times the elements of $\Cal F_3$ and finally differentiating five times the elements
of $\Cal F_5$:
$$
z^{l_1,1}_1=X_1z^{l_1,1},\ z^{1,3}_{r_3+1}=X_1z^{1,3}_{r_3},\ z^{l_5,5}_{r_5+1}=X_1z^{l_5,5}_{r_5}
$$   
where
$$
\aligned
1\leq l_1\leq 2,\ 1\leq l_5\leq 2,\cr
0\leq s_3\leq 2,\ 0\leq s_5\leq 4
\endaligned
$$
and $z^{l_j,j}_0:=z^{l_j,j}$. 
We thereby obtain the functions
{\eightpoint
$$
\aligned
&x=e^{x_1},z^{1,1}=x_3,z^{2,1}=x_{18},z^{1,1}_1=x_5-x_{21}-x_{20}+x_4-x_8-e^{x_1},z^{2,1}_1=x_{20},\cr
&z^{1,3}=x_6,z^{1,3}_1=-2x_{12}+x_7+2,z^{1,3}_2=x_8-x_{13}+x_{14}+e^{x_1},\cr
&\hskip 150 pt z^{1,3}_3=x_9+2(x_{20}+x_{21}-2x_4+2x_8+x_{16})+3e^{x_1},\cr
&z^{1,5}=x_{10},z^{2,5}=-x_2+x_{11},z^{1,5}_1=1-(x_{12}+x_{13}+x_{14}),z^{2,5}_1=1/2(x_{13}+x_{14}),\cr
&z^{1,5}_2=1/2(x_{14}-x_{13}-4x_{15}),z^{2,5}_2=x_{15},z^{1,5}_3=x_{16}-2x_{17},z^{2,5}_3=x_{17},\cr
&z^{1,5}_4=x_{18}-2(x_2+x_{19}+x_8+e^{x_1}),z^{2,5}_4=x_{19}-x_8-e^{x_1},\cr
&z^{1,5}_5=2(x_4-x_8-e^{x_1}-x_{21})-x_{20}, z^{2,5}_5=2(-x_4+x_8+e^{x_1})+x_{20}+x_{21},
\endaligned
$$}
\noindent in accordance with {\tt Contact A}. By Theorem 4.2, these functions define a local diffeomorphism 
$\psi\:\Bbb R^{21}\to J^{\langle 2,0,1,0,2\rangle}$ satisfying
$$
\psi_*\CV=\Cal C\langle 2,0,1,0,2\rangle
$$
where $J^{\langle 2,0,1,0,2\rangle}$ denotes the partial prolongation of $J^1(\Bbb R,\Bbb R^5)$ in which two variables remain at {\it order 1}, one variable is prolonged to {\it order 3}, two are prolonged to {\it order 5}. Finally, the contact distribution on $J^{\langle 2,0,1,0,2\rangle}$ has the form
$$
\aligned
&\Cal C\langle 2,0,1,0,2\rangle=\{\P x+\sum_{l_1=1}^2 z^{l_1,1}_1\P {z^{l_1,1}}+
\sum_{h_3=0}^2 z^{1,3}_{h_3+1}\P {z^{1,3}_{h_3}}+\sum_{l_5=1}^2\sum_{h_5=0}^4 z^{l_5,5}_{h_5+1}\P {z^{l_5,5}_{h_5}},\cr
& \hskip 250 pt \P {z^{1,1}_1}, \P {z^{2,1}_1}, \P {z^{1,3}_3}, \P {z^{1,5}_5}, \P {z^{2,5}_5}   \},
\endaligned
$$
using the notational conventions in the proof of Theorem 4.2.

\newpage

{\it Example 5.2.} [{\bf Nonlinear Control Theory}]. In geometric control theory one considers a car moving in the $xy$-plane (see [8]) modelled as follows. The {\it state} of the car is described by four variables $(x,y,\theta,\phi)$. The ordered pair $(x,y)$ gives the coordinates on the $xy$-plane of the centre of the rear axle. The variable $\theta$ is the angle between the $x$-axis fixed on the plane and the vertical axis $V$ of the car, running perpendicular to the axles. The variable $\phi$ is the angle the front wheels make relative to $V$.



Assuming the wheels do not slip as the car moves in the plane, then one obtains the Pfaffian system
$$
-\sin\theta ~dx +\cos\theta ~dy=0,\ L\cos\phi ~d\theta-\sin\phi~(\cos\theta~dx+\sin\theta~dy)=0.
$$    
This in turns leads to the control system $I=\{\o_1=0,\ldots,\o_4=0\}$, known as the {\it kinematic car}, where
$$
\o_1=dx-u^1\cos\theta~dt,\o_2=dy-u^1\sin\theta~dt,\o_3=d\theta-\frac{u^1}{L}\tan\phi~dt,\o_4=d\phi-u^2~dt,
$$
and $L$ is the length from the rear to the front axle. Here $u^1$ models the speed of the point $(x,y)$ and $u^2$ the speed at which the front wheels swivel. Variables $u^1,u^2$ are the {\it controls}, for prescribing these as functions of time 
$t$ gives a system of ordinary differential equations for the state of the car. 

Equivalently, control system $I$ defines the sub-bundle $\CK=I^\perp\subset T~(\Bbb R_t\times M)$, given by
$$
\CK=\{\P t + u^1(\cos\theta~\P x+\sin\theta~\P y +\frac{1}{L}\tan\phi~\P \theta)+u^2~\P \phi,\P {u^1},\P {u^2}\}
$$
where $M=\Bbb R^2_{(x,y)}\times \Bbb R^2_{(u_1,u_2)} \times S^1\times S^1$ is the manifold of states and controls. An important question in nonlinear control theory is: when can a nonlinear control system be ``linearised" by a {\it static feedback transformation} and more generally a {\it dynamic feedback transformation}? We will study these questions for the example of the kinematic car as an application of our general algorithm {\tt Contact}.

We begin by showing that $\CK$ is, in fact, a Goursat bundle of type $\tau=\langle 1,0,1\rangle$. 
The refined derived type of $\CK$ is
$$
[[3,0],[5,2,3],[6,4,4],[7,7]].
$$
Consequently the derived length is 3 and $\text{\rm deccel}(\CK)=\langle 1,0,1\rangle$. The filtration induced by $\CK$ is
$$
\ch \CK^{(1)}_0\subset\ch \CK^{(1)}\subset\ch \CK^{(2)}\subset T\Bbb R^7
$$
where
$$
\ch\CK^{(1)}_0=\{\P {u_1},\P {u_2}\},\ \ch\CK^{(1)}=\{\P {u_1},\P {u_2},\P t\},\ 
\ch\CK^{(2)}=\{\P {u_1},\P {u_2},\P t,\P {\phi}\}.\eqno(5.4)
$$
Since  the distributions in (5.4) are all integrable and $\r_3=1$, we conclude that $\CK$ is a Goursat bundle of type
$\langle 1,0,1\rangle$. By Theorem 4.1, $\CK$ is locally equivalent to the contact distribution 
$\Cal C\langle 1,0,1\rangle$. This settles the recognition question for $\CK$. We go on to find an equivalence using {\tt Contact}. 
 
For this, according to {\tt Contact B}, we compute distribution $\Pi^3$. 
The invariants of $\ch\CK^{(2)}$ are $x,y,\theta$. Any one of these may be taken to be the independent variable. 
If we choose $x$ for this purpose then we  take
$$
X=\P x+\frac{1}{u^1\cos\theta}(\P t+u^1\sin\theta~\P y+L^{-1}u^1\tan\phi~\P \theta+u^2~\P \phi)
$$  
for the total differential operator. By (4.13), we find that 
$$
\Pi^3=\{\P {u_1},\P {u_2},\P t,\P \phi,\P \theta\}.
$$
The filtration of $T^*\Bbb R^7$ induced by $\Cal K$ is therefore
$$
{\Pi^3}^\perp\subset\Xi^{(2)}\subset\Xi^{(1)}\subset\Xi^{(1)}_0\subset T^*\Bbb R^7
$$
where
$$
{\Pi^3}^\perp=\{dx,dy\},\ \Xi^{(2)}=\{dx,dy,d\theta\},\ \Xi^{(1)}=\{dx,dy,d\theta,d\phi\},\ 
\Xi^{(1)}_0=\{dx,dy,d\theta,d\phi,dt\}.
$$
It follows that the fundamental bundles are
$\O_1=\{dt\}$ and $\O_3=\{dy\}$
and the fundamental functions are therefore
$
\Cal F_1=\{t\},\ \Cal F_3=\{y\}.
$
Consequently, by {\tt Contact B} the map $\psi\:\Bbb R^7\to J^{\langle 1,0,1\rangle}$ defined by
$$
x,z^{1,1}=t,z^{1,1}_1=Xt,z^{1,3}=y,z^{1,3}_1=Xy,z^{1,3}_2=X^2y,z^{1,3}_3=X^3y\eqno(5.5)
$$
is an equivalence, according to Theorem 4.3. That is, $\psi_*\CK=\Cal C\langle 1,0,1\rangle$.

While (5.5) is certainly an identification with a partial prolongation, it is not of much use in control theory. This is because it is not a feedback equivalence. That is, (5.5) does not respect the special role played by the time coordinate 
$t$ in the original control system, nor the distinction between the roles played by the state and control variables. We  pause briefly to describe the class of transformations that preserves the set of all control systems.

Let 
$$
\frac{d{\bold x}}{dt}=f(t,{\bold x},{\bold u})\eqno(5.6)
$$
be a control system, where ${\bold x}$ denotes the state variables and ${\bold u}$ the control variables.
A local diffeomorphism of the form
$$
t\mapsto t,\ {\bold x}\mapsto \Phi({\bold x}),\ {\bold u}\mapsto \Psi({\bold x},{\bold u})\eqno(5.7)
$$
is said to be a {\it (static) feedback transformation}. A question of interest is: does there exist a static feedback transformation that identifies a given control system with some partial prolongation of $\Cal C^{(1)}_q$ for some $q$
or, as it is more commonly known in the control theory literature, a Brunovsky normal form?  

Writing, (5.5) out explicitly shows that it is {\it not} a static feedback equivalence. In fact, an elegant result of Sluis [5], applied to the kinematic car system $I$, proves that no feedback equivalence exists. Sluis' theorem relies on the GS algorithm and consequently applies to time-invariant control systems such as the kinematic car. In fact, a 
stronger necessary condition for static feedback equivalence can be derived via Theorems 4.2 and 4.3. 

\vskip 5 pt

{\smc Theorem 5.1.} {\it Let $\CK=\{\P t+f(t,{\bold x},{\bold u})\P {\bold x},\P {\bold u}\}\subset 
T~(\Bbb R_t\times M)$ arise from any smooth control system (5.6), where $M$ is the manifold of states and controls.   
A necessary condition in order that $\CK$ be static feedback equivalent to a Brunovsky normal form is that 
\item{\rm [i]} $\CK$ is a Goursat bundle 
\item{\rm [ii]} If $k$ is the derived length of $\CK$ then $dt\in{\ch\CK^{(k-1)}}^\perp$ if $\r_k=1$ or 
$dt\in\Upsilon_{\widehat{\S}_{k-1}}(\CK^{(k-1)})$ if $\r_k>1$.}

\vskip 3 pt
{\it Proof.} Suppose there is static feedback equivalence $\vartheta$ identifying $\CK$ with some Brunovsky normal form. Every such normal form is a partial prolongation $\Cal C\langle \tau\rangle$ of $\Cal C^{(1)}_q$, where 
$\tau=\langle \r_1,\r_2,\ldots,\r_k\rangle$, $\r_k\geq 1$. In the case $\r_k>1$, the independent variable $x$ of $\Cal C\langle \tau\rangle$ is an invariant of its resolvent bundle. Consequently, 
$\vartheta^*x=t$ is an invariant of the resolvent bundle $\Cal R_{\S_{k-1}}(\CK^{(k-1)})$ determined by $\CK$.
In the case $\r_k=1$, $x$ is an invariant of the Cauchy bundle $\ch\Cal C\langle \tau\rangle^{(k-1)}$ of the ${k-1}^{th}$ derived bundle $\Cal C\langle\tau\rangle^{(k-1)}$. Hence, $\vartheta^*x=t$ is an invariant of $\ch\CK^{(k-1)}$. \hfill\hfill \qed

\vskip 5 pt
While for the kinematic car, $\CK$ is certainly a Goursat bundle, it does not satisfy condition [ii] of Theorem 5.1. We deduce that $\CK$ is not static feedback linearisable, in agreement with Sluis' theorem.

However, as is well known for this example, a certain {\it Cartan prolongation} of $\CK$ {\it is} static feedback linearisable. We will not give precise definitions for this, referring the interested reader to [5] for the details. We merely want to present a Cartan prolongation of $\CK$ and show how to apply algorithm {\tt Contact}
to determine a static feedback linearisation of the Cartan prolonged distribution. 
  
We obtain a Cartan prolongation of the kinematic car system as follows. Define a new control system 
$$
\frac{dx}{dt}=u^1\cos\theta,\frac{dy}{dt}=u^1\sin\theta,\frac{d\theta}{dt}=\frac{u^1}{L}\tan\phi,
\frac{d\phi}{dt}=u^2, \frac{du^1}{dt}=w^2,\frac{dw^2}{dt}=v_1\eqno(5.8)
$$
by ``twice differentiating $u^1$". In control system (5.8) the coordinate $u^1$ has become a state variable and the new control variables are $v^1$ and $u^2$. To see the significance for control theory of this admittedly {\it ad hoc} construction one needs to check for the existence of a static feedback equivalence for system (5.8). We begin by changing notation slightly and examining the sub-bundle $\text{\rm pr}~\CK\subset T~(\Bbb R_t\times{\bar M})$, defined by (5.8), where $\bar M$ is the manifold of new states and controls. Setting $u^1=w^1$ and $u^2=v^2$, we have
$$
\text{\rm pr}~\CK=\{T=\P t+w^1(\cos\theta~\P x+\sin\theta~\P y+\frac{\tan\phi}{L}~\P \theta)+
v^2~\P \phi+w^2\P {w^1}+v^1\P {w^2},\P {v^1},\P {v^2}\}.
$$
The refined derived type is
$$
[[3,0],[5,2,2],[7,4,4],[9,9]].
$$
Proposition 3.1 shows that this is the derived type of the partial prolongation $\Cal C\langle 0,0,2\rangle$, that is, the  total prolongation $\Cal C^{(3)}_2$. A calculation reveals that 
$\text{\rm pr}~\widehat{\CK}^{(2)}:= \text{\rm pr}~\CK^{(2)}/\ch\text{\rm pr}~\CK^{(2)}$ is spanned by
$$
\text{\rm pr}~\widehat{\CK}^{(2)}=\widehat{\pi}~\{~{\G}=\P t+w^1(\cos\theta~\P x+\sin\theta~\P y),\P {w^1},\P \theta\}
$$
and that the polar matrix of a line $E=[a^1\widehat{\pi}({\G})+a^2\widehat{\pi}(\P {w^1})+a^3\widehat{\pi}(\P \theta)]$ is
$$
\left(\matrix -a^2 &a^1 & 0\\
              -a^3 & 0 & a^1\endmatrix\right)
$$
with respect to the non-zero structure
$$
\widehat{\d}\big(\widehat{\pi}(\P {z^1}),\widehat{\pi}({\G})\big)=
\big\langle\widehat{\pi}(\cos\theta~\P x+\sin\theta~\P y)\big\rangle,\ 
\widehat{\d}\big(\widehat{\pi}({\G}),\widehat{\pi}(\P {\theta})\big)=
\big\langle \widehat{\pi}\big(z^1(\sin\theta~\P x-\cos\theta~\P y)\big)\big\rangle.
$$
Here, again, $\widehat{\d}$ is the structure tensor of $\text{\rm pr}~\widehat{\CK}^{(2)}$,  
$$
\widehat{\pi} : T(\Bbb R_t\times\bar{M})\to T(\Bbb R_t\times\bar{M})/\ch\text{\rm pr}~\CK^{(2)}
=:\widehat{T}(\Bbb R_t\times\bar{M})
$$
is the natural projection and for $Z\in T(\Bbb R_t\times\bar{M})$, $\big\langle\widehat{\pi}(Z)\big\rangle$ denotes the 
element of the quotient bundle
$$
\widehat{T}(\Bbb R_t\times\bar{M})/\text{\rm pr}~\widehat{\CK}^{(2)}
$$ 
with representative $\widehat{\pi}(Z)$.
It follows that  $\text{\rm Sing}(\text{\rm pr}~\widehat{\CK}^{(2)})=
\Bbb P\{\widehat{\pi}(\P {w^1}),\widehat{\pi}(\P {\theta})\}$ and so, 
the resolvent bundle of $\text{\rm pr}~\CK^{(2)}$ is
$$
\Cal R_{\widehat{\S}_{2}}(\text{\rm pr}~\CK^{(2)})=\ch\CK^{(2)}\oplus\{\P {w^1},\P {\theta}\},
$$
which is integrable. The above data shows that $\text{\rm pr}~\CK$ is a Goursat bundle of type $\langle 0,0,2\rangle$.
By Theorem 4.1, there is an equivalence identifying $\text{\rm pr}~\CK$ with the contact distribution 
$\Cal C\langle 0,0,2\rangle$, that is, a Brunovsky normal form. But is it a static feedback equivalence?
To find out, we use algorithm {\tt Contact}. 

The only non-empty fundamental bundle in this case is the one of order 3
$$
\Upsilon_{\widehat{\S}_{2}}(\text{\rm pr}~\CK^{(2)})={\Cal R_{\widehat{\S}_{2}}(\text{\rm pr}~\CK^{(2)})}^\perp
=\{dt,dx,dy\}.
$$ 
So in this case, $dt\in\Upsilon_{\widehat{\S}_{2}}(\text{\rm pr}~\CK^{(2)})$ and we may choose $z^{1,3}=x,z^{2,3}=y$
and
$$
z^{1,3}_1=Tz^{1,3},z^{2,3}_1=Tz^{2,3},z^{1,3}_2=Tz^{1,3}_1,z^{2,3}_2=Tz^{2,3}_1,z^{1,3}_3=Tz^{1,3}_2,z^{2,3}_3=Tz^{2,3}_2.
$$
It is a simple matter to verify that these are indeed contact coordinates as predicted by Theorem 4.2. Moreover, an inspection of the formulas 
{\eightpoint
$$
\aligned
&z^{1,3}_1=w^1\cos\theta,\ z^{2,3}_1=w^1\sin\theta,\cr
&z^{1,3}_2= w^2\cos\theta-\frac{(w^1)^2}{L}\sin\theta\tan\phi,\ 
z^{2,3}_2=w^2\sin\theta+\frac{(w^1)^2}{L}\cos\theta\tan\phi,\cr 
&z^{1,3}_3=\frac{-(w^1)^3\cos\theta\sin^2\phi-3w^1w^2L\sin\phi\cos\phi\sin\theta-(w^1)^2v^2L\sin\theta+
v^1L^2\cos\theta\cos^2\phi}{L^2\cos^2\phi}\cr
&z^{2,3}_3=\frac{-(w^1)^3\sin\theta\sin^2\phi+3w^1w^2L\sin\phi\cos\phi\cos\theta+(w^1)^2v^2L\cos\theta+
v^1L^2\sin\theta\cos^2\phi}{L^2\cos^2\phi}
\endaligned
$$
}
shows that they define a static feedback equivalence. Because of the projection $\bar{M}\to M$, any integral manifold of
$\text{\rm pr}~\CK$ maps to a unique integral manifold of $\CK$. This, together with the fact that there is a static feedback equivalence for $\text{\rm pr}~\CK$ implies that the kinematic car example is {\it dynamic} feedback linearisable
in the language of nonlinear control theory. It is an important and largely open problem to geometrically characterise the class of nonlinear control systems which are dynamic feedback linearisable.

\vskip 20 pt
{\eightpoint

\centerline{\it References}

\item{[1]} R.L. Bryant, {\it Some Aspects of the Local and Global Theory of Pfaffian Systems}, Ph.D thesis, University of North Carolina, Chapel Hill, 1979
\item{[2]} R.L. Bryant, S-S. Chern, R.B. Gardner, H. Goldschmidt, P. Griffiths, {\it  Exterior Differential Systems},
MSRI Publications, \allowbreak Springer-Verlag, New York, 1991
{
\item{[3]} R.B. Gardner, W.F. Shadwick, The GS algorithm for exact linearization to Brunovsk\'y normal form, {\it IEEE Trans. Automat. Control}, 37(2) (1992), 224-230

\item{[4]} W. Respondek, W. Pasillas-Lepine,
Contact systems and corank 1 involutive subdisributions, {\it Acta Appl. Math.}, {\bf 69} (2001), 105-128
\item{[5]} W. Sluis, {\it Absolute Equivalence and its Applications to Control Theory}, Ph.D Thesis, University of Waterloo, Waterloo, Canada, 1992  
\item{[6]} P. Vassiliou, A constructive generalised Goursat normal form, {\tt  Mathematics Archive}, {\tt math.DG/0404377}, April 2004
\item{[7]} E. von Weber, Zur invariantentheorie der systeme Pfaff'scher gleichungen, \newline {\it Berichte Verhandlungen
der Koniglich Sachsischen Gesellshaft der Wissen-\newline shaften Mathematisch Physikalische Klasse, Leipzig}, 50: 207-229, (1898)
\item{[8]} G. Wilkens, Finsler geometry in low dimensional control theory, in D. Bao, S-S. Chern, Z. Shen (Eds), Finsler Geometry, 
{\it Contemporary Mathematics Vol. 196}, American Mathematical Society, 1996

\item{[9]} K. Yamaguchi, Contact geometry of higher order, {\it Japan J.} 
 {\it Math.}, 8, (1982), 109-176; Geometrization of jet bundles, {\it Hokkaido Math. J.}, 12, (1983), 27-40

\vskip 10pt
{\bf Acknowledgments}: The computations in this paper were immensely assisted by use of the general purpose differential geometry package {\tt VESSIOT} developed at Utah State University by Ian M. Anderson and colleagues.    
\vskip 10 pt
}

\enddocument
\bye